\documentclass[12pt]{article}
\usepackage{amssymb,theorem,a4wide}
\newcommand{\X}{{\cal X}}
\newcommand{\Xcoarse}{\overline{\X}}
\newcommand{\Xcan}{{\Xcoarse^{\mathrm{an}}}}
\newcommand{\Xpc}{Q}
\newcommand{\M}{{\cal M}}
\newcommand{\etabar}{{\overline{\eta}}}
\newcommand{\Spec}{\mathop{\mathrm{Spec}}}
\newcommand{\shfO}{{\cal O}}

\newcommand{\NN}{{\mathbb{N}}}
\newcommand{\FF}[1]{{\mathbb{F}_{#1}}}
\newcommand{\FFpbar}{{\overline{\mathbb{F}}_p}}
\newcommand{\ZZ}{{\mathbb{Z}}}
\newcommand{\ZZp}{{\ZZ_p}}
\newcommand{\ZZl}{{\ZZ_l}}
\newcommand{\ZZpnr}{\mathbb{Z}_p^{\mathrm{nr}}}
\newcommand{\ZZhat}{\hat\ZZ}
\newcommand{\QQ}{{\mathbb{Q}}}
\newcommand{\QQbar}{{\overline{\QQ}}}
\newcommand{\QQp}{{\QQ_p}}
\newcommand{\QQpnr}{\mathbb{Q}_p^{\mathrm{nr}}}
\newcommand{\QQpbar}{{\overline{\QQ}_p}}
\newcommand{\QQl}{{\QQ_l}}
\newcommand{\QQln}[1]{\QQ_l^{#1}}
\newcommand{\QQlbar}{{\overline{\QQ}_l}}
\newcommand{\CC}{{\mathbb{C}}}

\newcommand{\Lo}{o\hskip1pt}

\newcommand{\rH}{\mathrm{H}}
\newcommand{\rR}{\mathrm{R}}
\newcommand{\et}{{\mathrm{\acute et}}}

\newcommand{\trace}{\mathrm{Tr}}

\newcommand{\Gal}{\mathrm{Gal}}
\newcommand{\Frob}{\mathrm{Frob}}

\newcommand{\Dst}[1]{\mathrm{D}_{\mathrm{st},#1}}
\newcommand{\ssrep}[1]{\underline{\mathrm{Rep}}_{\mathrm{st}}(#1)}
\newcommand{\unrrep}[1]{\underline{\mathrm{Rep}}_{\mathrm{unr}}(#1)}
\newcommand{\MFNf}[1]{\underline{\mathrm{MF}}_{#1}^f(\varphi,N)}
\newcommand{\Fil}{\mathrm{Fil}}

\newcommand{\Hom}{\mathrm{Hom}}
\newcommand{\Aut}{\mathrm{Aut}}
\newcommand{\Ext}{\mathrm{Ext}}

\newcommand{\DR}{{\mathrm{DR}}}
\newcommand{\Ddr}{{\mathrm{D}_\DR}}

\newcommand{\Affine}{\mathbf{A}}
\newcommand{\MU}{\mu\kern-.5em\mu}

\newcommand{\map}[2]{\;\smash{\mathop{#1}\limits^{#2}}\;}

\theorembodyfont{\sl}
\newtheorem{theorem}[subsection]{Theorem}
\newtheorem{proposition}[subsection]{Proposition}
\newtheorem{lemma}[subsection]{Lemma}
\newtheorem{corollary}[subsection]{Corollary}

\newenvironment{assumption}[1]{$$\hbox{(#1)}\quad
  \vtop\bgroup\hsize.8\hsize\noindent\strut\ignorespaces}{%
  \egroup$$\ignorespacesafterend}

\newenvironment{prf}[1]{\trivlist
\item[\hskip \labelsep{\it
#1.\hspace*{.3em}}]}{~\hspace{\fill}~$\square$\endtrivlist}
\newenvironment{proof}{\begin{prf}{\bf Proof}}{\end{prf}}

\begin{document}

\author{Theo van den Bogaart \and Bas Edixhoven\footnote{Both authors
were partially supported by the European Marie Curie Research Training
Network ``arithmetic algebraic geometry'', contract
MRTN-CT2003-504917. Address: Mathematisch Instituut, Universiteit
Leiden, Postbus 9512, 2300\ RA\ \ Leiden, The Netherlands,
bogaart@math.leidenuniv.nl, edix@math.leidenuniv.nl}} 

\title{Algebraic stacks whose number of points over finite fields is a
polynomial} \date{November 1, 2008} 

\maketitle

\section{Introduction}
The aim of this article is to investigate the cohomology ($l$-adic as
well as Betti) of schemes, and more generally of certain algebraic
stacks, $\X$, say, that are proper and smooth over the spectrum of~$\ZZ[1/n]$ 
for some $n\ge1$ and have the
property that there exists a polynomial $P$ with coefficients in~$\QQ$
such that for every finite field $\FF{q}$ of characteristic not dividing~$n$ 
we have $\#\X(\FF{q}) =
P(q)$. For the precise definitions and conditions the reader is
invited to read the rest of the article, at least up to the statement
of Theorem~\ref{thm:main}. Under those conditions, we prove that for
all prime numbers~$l$ the \'etale cohomology
$\rH(\X_{\QQbar,\et},\QQl)$, considered as a representation of the
absolute Galois group of~$\QQ$, is as expected: zero in odd degrees,
and, after semi-simplification, a direct sum of $\QQl(-i)$ in degree~$2i$, 
with the number of
terms equal to the coefficient $P_i$ of~$P$.
For $n=1$ we prove that $\rH(\X_{\QQbar,\et},\QQl)$ is semi-simple.
Our main tools here are
Behrend's Lefschetz trace formula in~\cite{Behrend} and $l$-adic Hodge
theory combined with the fact that $\ZZ$ has no nontrivial unramified
extensions. Finally, using comparison theorems from $l$-adic Hodge
theory, we obtain a corollary which says that, under the extra
assumption that the coarse moduli space of~$\X$ is a quotient by a
finite group, the Betti cohomology $\rH(\X(\CC),\QQ)$ with its Hodge
structure is as expected: zero in odd degree, and $\QQ(-i)^{P_i}$ in
degree~$2i$.

The results in this article are motivated by a question by Carel Faber
on potential applications to some moduli stacks $\M_{g,n}$ of stable
$n$-pointed curves of genus~$g$. These stacks are proper and smooth
over~$\ZZ$, and they do also satisfy the extra hypotheses of the
corollary by results of Pikaart and Boggi~\cite{Pikaart-Boggi}.  We
are told that $\#\M_{g,n}(\FF{q})$ is a polynomial in~$q$ for all
pairs of the form $(0,n)$ with $n\geq3$, $(1,n)$ with $1\leq n\leq
10$, $(2,n)$ with $0\leq n\leq 5$ (probably even up to~$n=9$), and
$(3,n)$ with $0\leq n\leq 3$ (and probably more).  For genus 2~and~3
these results are due to Jonas Bergstr\"om.

As this article is motivated by its application to certain~$\M_{g,n}$,
we have not made an effort to make our results as general as
possible. In particular, we have not tried to generalise comparison
theorems from $l$-adic Hodge theory from schemes to stacks.

We hope that this article will be of help to those computing the
rational Hodge structure on the cohomology of
certain~$\M_{g,n}$. Counting points, using a suitable stratification,
could be easier than having to compute the cohomology, using the same
stratification. We apologise for our lack of expertise in the fields
of algebraic stacks and $l$-adic Hodge theory.  Readers with more
competence in these areas will probably find the contents of this
article rather straightforward and the proofs too elaborate. But there
seems to be a lack of `well-known facts' in the literature and we have
tried hard to give precise references and proofs understandable also
to the non-expert.

\subsection*{Acknowledgements}
We thank Carel Faber for asking us his question, and the organisers of
the 2004~Texel Island conference where this happened. We also thank
Laurent Moret-Bailly for giving us the opportunity to ask him some
questions about stacks and Gerard van der Geer for his comments.

\subsection*{Terminology, conventions}
Concerning stacks, our terminology is that of~\cite{LMB}. In particular,
algebraic stacks are by definition quasi-separated.

Let $k$ be a finite field. If $\X$ is a Deligne-Mumford stack 
of finite type over~$\ZZ$,
define its \emph{number of points over~$k$} to be
$$\#\X(k)=\sum_\xi\frac{1}{|\Aut(\xi)|}\;,$$
where the sum is over representatives of isomorphism classes of objects 
in~$\X(k)$. Here $\Aut(\xi)$ denotes the finite group of automorphisms 
of~$\xi$.

If $G$ is a topological group, by a $G$-representation (over~$\QQl$)
we shall mean a continuous representation of $G$ on a finite
dimensional $\QQl$-vector space equipped with the $l$-adic
topology. We use the same notation for a representation and its
underlying vector space.  For any~$G$ and $n\ge0$, the symbol
$\QQln{n}$ denotes the trivial $n$-dimensional $G$-representation.

\subsection*{Updates}
After version~1 of this text appeared on arxiv, some related results
that use the relations between point counting over finite fields,
\'etale cohomology and Hodge numbers have been brought to our
attention. In~\cite{Ito}, Ito uses these tools to reprove that the
stringy $E$-function does not depend on the choice of a resolution of
singularities. In \cite{KL1} Kisin and Lehrer use point counting over
finite fields to determine the character of a finite group acting on
the cohomology of a variety. In~\cite{KL2} they study varieties, not
necessarily proper or smooth, whose number of points over finite
fields is given by a polynomial, as in our case, and give two
sufficient conditions for the cohomology to be of Tate type. A result
like our Lemma~\ref{lemma:semisimple}, and its application to global
Galois representations, is in the article~\cite{KW} by Kisin and
Wortmann. In~\cite{BT} Bergstr\"om and Tommasi have adapted our proofs
to the equivariant case for the action by a finite group; their
article shows that $\M_4$ and all strata of its boundary satisfy the
polynomiality condition.  Bergstr\"om has also written two other
articles (\cite{BergA} and~\cite{BergB}) involving the polynomiality
condition for various moduli spaces of curves.

The first author has proved, in~\cite[Cor.~8.12]{Bogaart}, the de Rham
comparison theorem from $p$-adic Hodge theory for Deligne-Mumford
stacks that are smooth and proper over complete discrete valuation
fields with perfect residue field of characteristic~$p$. This enables
us to remove the extra condition that the coarse moduli space is a
certain quotient space in Corollary~\ref{cor:main} of the last section
of this article (see the footnote there for more detail).

We thank Moret-Bailly for bringing to our attention the fact that our
proof gives a more general result than the one in the statement of the
original theorem. Part of the conclusions of the main theorem are
valid under the assumption that the stack is smooth and proper over a
dense open part of~$\Spec\ZZ$; originally, we had the stronger
assumption that it had to be smooth and proper over the whole
of~$\Spec\ZZ$.

\section{Results}

\begin{theorem}\label{thm:main}
Let $\X$ be a Deligne-Mumford stack over~$\ZZ$.
Let $d\ge0$ and assume that $\X$~is proper, smooth
and of pure relative dimension~$d$ over some non-empty open subscheme~$U$
of~$\Spec\ZZ$.
Let $S$ be a set of primes of Dirichlet density~1. 
Assume~---
\begin{assumption}{$*$}
there exists a polynomial $P(t)=\sum_{i\ge0}P_it^i$, with $P_i\in\QQ$, such 
that
$$\#\X(\FF{p^n})=P(p^n)+\Lo(p^{nd/2})\qquad(n\to\infty)$$
for all $p\in S$.
\end{assumption}
Then the degree of $P(t)$ is~$d$, and there exists a unique such 
polynomial satisfying $P_i=P_{d-i}$ for all $0\le i\le d$.
Suppose $P(t)$ is of this form.
Then it has non-negative integer coefficients and satisfies
$\#\X(\FF{p^n})=P(p^n)$ for all primes~$p$ in~$U$ and all~$n\ge1$.
Furthermore, for all primes~$l$ and all $i\ge0$ there is an isomorphism
of $\Gal(\QQbar/\QQ)$-representations
$$\rH^i(\X_{\QQbar,\et},\QQl)^{\rm ss}\simeq
  \cases{0,&if $i$ is odd;\cr
  \QQl(-i/2)^{P_{i/2}},&if $i$ is even.}$$
Here $\rH^i(\X_{\QQbar,\et},\QQl)^{\rm ss}$ denotes the semi-simplification
of $\rH^i(\X_{\QQbar,\et},\QQl)$.
If $U=\Spec\ZZ$, then $\rH^i(\X_{\QQbar,\et},\QQl)^{\rm ss}\simeq
\rH^i(\X_{\QQbar,\et},\QQl)$.
\end{theorem}

We remark that part of the theorem can also be stated in terms of the
coarse moduli space associated to~$\X$. Indeed, the number of points
over a finite field of a Deligne-Mumford stack equals the number of
points of its coarse moduli space; and furthermore, the cohomologies
of both spaces with coefficients in a $\QQ$-algebra are the same.

\section{Some Results on Stacks}

Let us make two technical remarks.

In \cite[\S18]{LMB} the theory of constructible sheaves of
$\ZZ/{l^n\ZZ}$-modules over the smooth-\'etale site of an algebraic
$S$-stack is developed, where $S$ is a scheme. There is a
straightforward extension of this theory to constructible $l$-adic
sheaves, e.g., by working with projective systems of
$\ZZ/{l^n\ZZ}$-modules modulo torsion in the usual way. We will use
this without further comments.

Associated to a Deligne-Mumford stack~$\X$ are its \'etale topos
(denoted~$\X_\et$) and its smooth-\'etale topos. They however give the
same cohomology theory of constructible sheaves (see~\cite[\S12]{LMB},
especially Prop.~12.10.1). This justifies the fact that we will only
work with the \'etale topos of a Deligne-Mumford stack, but freely
cite results stated in terms of the other topos.

\medskip

The scheme counterpart of the following proposition is classical. By
lack of a precise reference, we have included a proof for the case of
stacks.

\begin{proposition}\label{prop:mon}
Let $\X$ be a Deligne-Mumford stack which is smooth and proper
over~$\ZZp$.  For every prime $l\not=p$ and every $i\ge0$, the
canonical map of $\Gal(\QQpbar/\QQp)$-representations
\begin{equation}\label{eq:spec}
\rH^i(\X_{\FFpbar,\et},\QQl)\to\rH^i(\X_{\QQpbar,\et},\QQl)
\end{equation}
is an isomorphism. In particular, $\rH^i(\X_{\QQpbar,\et},\QQl)$ is
unramified.
\end{proposition}

\begin{proof}
Denote by $\QQpnr$ the maximal unramified extension of $\QQp$
in~$\QQpbar$ and let $\ZZpnr$ be its ring of integers.  Set
$S=\Spec(\ZZpnr)$ and denote by $s$ resp.~$\eta$ its closed
resp. generic point. Let $\etabar\to\eta$ correspond to
$\QQpnr\to\QQpbar$.  Consider the natural morphisms
$$\X_\etabar \map{\rightarrow}{j} \X_S \map{\hookleftarrow}{i} \X_s\;.
$$ These maps induce continuous morphisms between the associated
\'etale sites.

Let $(U,u)$ be an \'etale neighbourhood of~$\X_S$ and let $j_U$ be the
pull-back of $j$ along~$u$. By \cite[18.2.1(i)]{LMB}, for every $q$ we
have $(\rR^qj_*\QQl)_{U,u}=\rR^q(j_U)_*\QQl$. As $U$ is smooth, it
follows that $j_*\QQl=\QQl$ and $\rR^qj_*\QQl=0$ if~$q\not=0$.  Hence
the Leray spectral sequence gives an isomorphism
$\rH^i(\X_S,\QQl)\to\rH^i(\X_\etabar,\QQl)$.  But on the other hand,
$\rH^i(\X_S,\QQl)$ is naturally isomorphic to $\rH^i(\X_s,\QQl)$; this
follows from the proper base change theorem (\cite[18.5.1]{LMB}) for $\X_S$
over~$S$ and the fact that $S$ is strictly local.
\end{proof}

\bigskip

The next topic is Poincar\'e duality for the $l$-adic cohomology of
certain stacks (see Prop.~\ref{prop:PD} below). We will obtain this by
considering the cohomologies of their associated coarse moduli spaces.

Let $\X$ be a separated Deligne-Mumford stack of finite type over an
algebraically closed field of characteristic zero.  We will denote by
$\Xcoarse$ its coarse moduli space and by $q\colon\X\to\Xcoarse$ the
corresponding mapping. Note that we can cover $\Xcoarse$ by \'etale
charts $U$ such that the pull-back of $U$ in $\X$ is the quotient
stack of an algebraic space by a finite
group~(\cite[Rem.~6.2.1]{LMB}).

\begin{lemma}\label{lemma:cohomcoarse}
For every $i$ the pull-back map
$$q^*\colon\rH^i(\Xcoarse_\et,\QQl)\longrightarrow\rH^i(\X_\et,\QQl)$$
is an isomorphism.
\end{lemma}

\begin{proof}
The lemma follows from the Leray spectral sequence once we have shown
that the canonical map $\QQl\to\rR q_*\QQl$ is an isomorphism.  This
question is \'etale local on~$\Xcoarse$ and therefore we may assume
that $\X=[V/G]$ for some algebraic space~$V$ equipped with an action
by a finite group~$G$. Denote by $p\colon V\to\X$ the canonical
morphism.  Note that $\QQl\simeq(p_*\QQl)^G$. As $p$ and $qp$ are
finite and $\QQ[G]$ is a semi-simple $\QQ$-algebra, we obtain
$$\rR q_*\QQl\simeq\rR q_*(p_*\QQl)^G\simeq((qp)_*\QQl)^G\simeq\QQl
\;.$$
\end{proof}

Now suppose that $\X$ is defined over~$\CC$ and smooth. Consider the
complex analytic space $\Xcan$ associated to~$\Xcoarse$. It can
naturally be equipped with the structure of a $V$-manifold, i.e.,
locally $\Xcan$ is the quotient of a connected manifold by a finite
group; c.f.,~\cite{Steenbrink}.

\begin{proposition}\label{prop:PD}\label{prop:correspondence}
Suppose $\X$ is a Deligne-Mumford stack which is smooth and proper
over~$\QQbar$ of pure dimension~$d$ for some~$d\ge0$.
\begin{itemize}
\item[i)] Suppose $\X$ is integral.  For an integer~$i$, consider the
cup product mapping
$$\rH^i(\X_\et,\QQl)\otimes_\QQl\rH^{2d-i}(\X_\et,\QQl)
  \longrightarrow\rH^{2d}(\X_\et,\QQl)\;.$$ Then the right-hand side
  is one-dimensional and the pairing thus obtained is perfect.
\item[ii)] Suppose $X$ is a smooth and proper $\QQbar$-scheme of pure
dimension~$d$ and let $f\colon X\to\X$ be a $\QQbar$-morphism which is
surjective and generically finite. Then for all~$i$, the induced map
$$f^*\colon\rH^i(\X_\et,\QQl)\longrightarrow \rH^i(X_\et,\QQl)$$ is
  injective.
\end{itemize}
\end{proposition}

\begin{proof}
By Lemma~\ref{lemma:cohomcoarse} and the comparison theorem between
Betti and \'etale cohomology, it suffices to show that in case~(i)
$$\rH^i(\Xcan,\QQ)\otimes\rH^{2d-i}(\Xcan,\QQ)\to
  \rH^{2d}(\Xcan,\QQ)$$ is a perfect pairing; and in case~(ii) that
$$\rH^i(\Xcoarse_\et,\QQl)\longrightarrow \rH^i(X_\et,\QQl)$$ is
  injective. Now the singular cohomology of a $V$-manifold satisfies
  Poincar\'e duality (\cite{Steenbrink}), from which these statements
  follow.
\end{proof}

\section{Proof of the Main Theorem}

We will now prove Thm.~\ref{thm:main}. In this section, all cohomology
is with respect to the \'etale sites. We begin with an analytic lemma
used in the course of the proof.

\begin{lemma}\label{lemma:analysis}
Let us be given the following integers: $d\ge0$, $d\le r\le 2d$, and
for $0\le i\le r$ also $d_i\ge0$. Furthermore, let $p>1$ be a real
number, let $P_i\in\QQ$ for all $i\ge0$, with $P_i=0$ for $i$ large,
and let $\alpha_{i,j}\in\CC$ for $0\le i\le r$ and $1\le j\le
d_i$. Assume $|\alpha_{i,j}|=p^{i/2}$ and
\begin{equation}\label{eq:estimate}
\sum_{i=0}^r (-1)^i\!\!\sum_{1\le j\le d_i}\alpha_{i,j}^n=
\sum_{i\ge0}P_ip^{ni} + \Lo(p^{nd/2})\qquad(n\to\infty)\;.
\end{equation}
Then $d_i=0$ for $i\ge d$ odd, $P_i=0$ for $i>r/2$, while for $d/2\le
i\le r/2$ we have $P_i=d_{2i}$; for these $i$ also
$\alpha_{2i,j}=p^i$.
\end{lemma}

\begin{proof}
The lemma follows by induction on~$r$. Indeed, assume that either
$r=d$ or that the lemma holds for $r-1$.  As $|\sum_{1\le j\le
d_i}\alpha_{i,j}^n|\le d_ip^{in/2}$, we have
$$\bigl|\sum_{i=0}^r(-1)^i\!\!\sum_{1\le j\le d_i}\alpha_{i,j}^n\bigr|
  =\bigl|\sum_{1\le j\le d_r}\alpha_{r,j}^n\bigr|+
  \Lo(p^{nr/2})\qquad(n\to\infty)$$ and also,
  using~(\ref{eq:estimate}), that $P_i=0$ for $i>r/2$.

Note that if $z$ is an element of a finite product~$(S^1)^s$ of
complex unit circles, then the closure of $\{z^n\mid n\ge1\}$ contains
the unit element. Hence for every $\epsilon>0$ there exists an
infinite subset $N\subset\NN$ such that for all $n\in N$ and for all
$i$~and~$j$ we have $|(\alpha_{i,j}p^{-i/2})^n-1| <\epsilon$ and in
particular
\begin{equation}\label{eq:Dir}
\bigl|\sum_{1\le j\le
d_r}(\alpha_{r,j}p^{-r/2})^n-d_r\bigr|<\epsilon'\;,
\end{equation}
with $\epsilon'=d_r\epsilon$.

Now first suppose $r$ is odd. Then
$$\bigl|\sum_{1\le i\le d_r}\alpha_{r,j}^n\bigr|=
  \Lo(p^{nr/2})\qquad(n\to\infty)\;,$$ which by~(\ref{eq:Dir}) implies
  $d_r=0$. Therefore, if $r=d$ we are done, while if $r>d$ we can
  apply the induction hypotheses.

So from now on suppose $r$ is even. From (\ref{eq:estimate}) follows
$$\bigl|\sum_{1\le j\le d_r}\alpha_{r,j}^n- P_{r/2}p^{n r/2}\bigr|=
  \Lo(p^{nr/2})\qquad(n\to\infty)\;,$$ or equivalently:
$$\lim_{n\to\infty} \bigl|\sum_{1\le j\le
  d_r}(\alpha_{r,j}p^{-r/2})^n-P_{r/2}\bigr|=0\;.$$ Together
  with~(\ref{eq:Dir}) this implies $d_r=P_{r/2}$.  In turn this easily
  leads to $\alpha_{r,j}=p^{r/2}$.

Now subtract $P_{r/2}p^{rn/2}=\sum_{1\le j\le d_r}\alpha_{r,j}^n$
from~(\ref{eq:estimate}) and if $r>d$ apply the induction hypotheses.
\end{proof}

Let $\X$, $d$, $U$ and $S$ be as in Thm.~\ref{thm:main} and let
$P(t)=\sum_{i\ge0}P_it^i$ be a polynomial satisfying~$(*)$. Without
loss of generality we assume that $P_i=P_{d-i}$ for all $0\le i\le d$
and that $S$~is contained in~$U$.  We also fix a prime~$l$.

\medskip

By \cite[Thm.~16.6]{LMB} and resolution of singularities, there exists
a smooth and proper $\QQ$-scheme $X$ of pure dimension~$d$ and a
surjective and generically finite map $f\colon X\to\X_\QQ$. By
removing a finite number of primes from $S$ if necessary, we may
assume that $X_\QQp$ extends to a smooth scheme over~$\ZZp$ for
all~$p\in S$.  As a consequence, for all $p\not=l$ in~$S$ the
representation $\rH^i(X_\QQpbar,\QQl)$ is unramified and we can
consider the action of Frobenius. By \cite{Deligne}, the eigenvalues
of Frobenius have complex absolute value~$p^{i/2}$. Using
Prop.~\ref{prop:correspondence}, we obtain the same conclusions for
the subrepresentation~$\rH^i(\X_{\QQpbar},\QQl)$.

Fix a prime $p\not=l$ in~$S$. By Behrend's Lefschetz trace formula
(see \cite{Behrend} or~\cite[Thm.~19.3.4]{LMB}),
\begin{equation}\label{eq:Lefschetz}
\sum_{i\ge0}(-1)^i\,\trace\bigl(\Frob^n,\rH^i(\X_\FFpbar,\QQl)\bigr)
=\#\X(\FF{p^n})
\end{equation}
for all~$n\ge1$. Let $\alpha_{i,1}$,~\ldots, $\alpha_{i,d_i}$ be the
complex roots of the characteristic polynomial of
Frobenius. Applying~$(*)$ and Prop.~\ref{prop:mon},
formula~(\ref{eq:Lefschetz}) becomes
$$\sum_{i=0}^{2d}(-1)^i\!\!\sum_{1\le j\le d_i}\alpha_{i,j}^n=
  P(p^n)+\Lo(p^{nd/2})\qquad(n\to\infty)\;.$$ From
  Lemma~\ref{lemma:analysis} we now obtain that $P(t)$ has degree~$d$
  and for all $d\le i\le 2d$ we have that for $i$ even $P_{i/2}=d_{i}$
  and $\alpha_{i,j}=p^{i/2}$, while $d_i=0$ for $i$ odd.  Using
  Poincar\'e duality (Prop.~\ref{prop:PD}) we obtain the same
  conclusions for all~$i$.  (Note that $P(t)$ is defined in such a way
  that $P_{i}=P_{d-i}$ for $0\le i\le d/2$.)

Hence $\rH^i(\X_\QQpbar,\QQl)$ vanishes for any odd~$i$, while for all
even~$i$ it has dimension $P_{i/2}$; in particular, the coefficients
of~$P(t)$ are non-negative integers. Furthermore,
\begin{equation}\label{eq:tr}
\trace\bigl(\Frob,\rH^{2i}(\X_\QQpbar,\QQl)\bigr)=
  \trace\bigl(\Frob,\QQl(-i)^{P_i}\bigr)\;.
\end{equation}
This holds for all primes~$p$ in the set $S\backslash\{l\}$.  But a
semi-simple, almost everywhere unramified,
$\Gal(\QQbar/\QQ)$-representation over~$\QQl$ is determined by the
trace of Frobenius on a set of primes of Dirichlet density~1 (a proof
of this is outlined in \cite[Prop.~2.6]{DDT}).  We conclude that the
semi-simplification of $\rH^{2i}(\X_\QQbar,\QQl)$ is isomorphic
to~$\QQl(-i)^{P_i}$.

By Prop.~\ref{prop:mon}, $\rH^{2i}(\X_\QQpbar,\QQl)$ is unramified for
every prime $p\not=l$ in~$U$. As (\ref{eq:tr}) then holds for all these
primes, Behrend's Lefschetz trace formula gives
$\#\X(\FF{p^n})=P(p^n)$. Changing $l$, we see that
this formula is valid for every prime~$p$ in~$U$.

\medskip

Now assume $U=\Spec\ZZ$.
All that remains to be proved is that $\rH^{2i}(\X_\QQbar,\QQl)$ is
semi-simple, or equivalently, that its $i$th Tate twist
$H=\rH^{2i}(\X_\QQbar,\QQl)(i)$ is semi-simple. Note that $H$ is
unramified outside~$l$ and that the semi-simplification of $H$ is
isomorphic to the trivial representation~$\QQln{P_i}$.

By \cite[16.6]{LMB} and \cite{DeJong}, there exists a finite extension
$K$ of~$\QQl$ inside~$\QQlbar$, a proper, semi-stable scheme $X$ of
pure relative dimension~$d$ over the ring of integers of~$K$ and a
surjective and generically finite $K$-morphism $f\colon X_K\to\X_K$.
By \cite[Thm.~1.1]{Tsuji}, $\rH^{2i}(X_\QQlbar,\QQl)$ is a semi-stable
representation of~$\Gal(\QQlbar/K)$. So it follows from
Prop.~\ref{prop:correspondence} that $\rH^{2i}(\X_\QQlbar,\QQl)$ and
hence also $H$ are potentially semi-stable.

\begin{lemma}\label{lemma:semisimple}
Let $n\ge1$ be an integer. Consider a short exact sequence of
$\Gal(\QQlbar/\QQl)$-representations
\begin{equation}\label{eq:sesrepr}
0\to\QQln{n}\to V\to\QQl\to0\;.
\end{equation}
If $V$ is potentially semi-stable, then $V$ is unramified.
\end{lemma}

\begin{proof}
By assumption, there is a finite extension $K$ of $\QQl$ inside
$\QQlbar$, such that the restriction of $V$ to $G_K:=\Gal(\QQlbar/K)$
is semi-stable. Fix such a~$K$ and denote by $K_0$ its maximal
unramified subfield relative to~$\QQl$. Denote by $\sigma$ the
automorphism of~$K_0$ obtained by lifting the automorphism $x\mapsto
x^l$ of the residue field of~$K_0$.

We will briefly recall some theory about semi-stable representations
and filtered~$(\varphi,N)$-modules; for more details, see~\cite{CF}.

Denote by $\ssrep{G_K}$ the category of semi-stable representations of
$G_K$ over~$\QQl$ and denote by $\unrrep{G_K}$ its full subcategory of
unramified representations. Let $\MFNf{K}$ be the category of (weakly)
admissible filtered $(\varphi,N)$-modules over~$K$. An object of
$\MFNf{K}$ is a finite dimensional $K_0$-vector space~$E$ equipped
with a $\sigma$-semi-linear bijection~$\varphi\colon E\to E$, a
nilpotent endomorphism $N$ of~$E$ and an exhaustive and separating
descending filtration~$\Fil^\cdot E_K$ on $E_K:=K\otimes_{K_0}E$. We
must have $N\varphi=l\varphi N$ and furthermore there is a certain
admissibility condition to be satisfied (c.f.,~\cite[\S3]{CF}).

All above categories are Tannakian (so in particular they are all
abelian $\QQl$-linear $\otimes$-categories). Fontaine has constructed
a functor~$\Dst{K}$ from $\ssrep{G_K}$ to $\MFNf{K}$ and the main
result of~\cite{CF} is that this is an equivalence of Tannakian
categories.

To the trivial one-dimensional representation~$\QQl$ corresponds
$\Dst{K}(\QQl)$, which is just $K_0$ equipped with the trivial maps
$\varphi=\sigma$, $N=0$ and filtration determined by $\Fil^0 K=K$
and~$\Fil^1 K=0$. By abuse of notation we will denote $\Dst{K}(\QQl)$
also by~$K_0$.

We obtain natural maps of $\QQl$-vector spaces
$$\Ext^1_{\unrrep{G_K}}(\QQl,\QQln{n})\map{\hookrightarrow}{i}
  \Ext^1_{\ssrep{G_K}}(\QQl,\QQln{n})\map{\rightarrow}{\sim}
  \Ext^1_{\MFNf{K}}(K_0,K_0^n)\;,$$ where the second map is the
  isomorphism induced by~$\Dst{K}$. We need to show that $i$ is an
  isomorphism. Since $\Ext^1$ is additive in the second variable, it
  suffices to treat the case~$n=1$. We will show that the dimension
  of $\Ext^1_{\unrrep{G_K}}(\QQl,\QQl)$ is at least as big as the
  dimension of the $\QQl$-vector space $\Ext^1_{\MFNf{K}}(K_0,K_0)$.

First we consider the extensions in $\MFNf{K}$. Let
\begin{equation}\label{eq:sesfpm}
0\to K_0 \to E \to K_0\to 0
\end{equation}
be a short exact sequence in~$\MFNf{K}$.  Choosing a splitting of the
short exact sequence of vector spaces underlying~(\ref{eq:sesfpm}), we
write $E=K_0\oplus K_0$ with $K_0\to E$ being $x\mapsto(x,0)$ and
$E\to K_0$ being the projection onto the second coordinate.  As the
induced sequence of filtered $K$-vector spaces is exact,
$\Fil^0E_K=E_K$ and $\Fil^1E_K=0$. Secondly, $N$ has the form
$\bigl({0\ \lambda\atop 0\ 0}\big)$ for some $\lambda\in K_0$. Then
$N\varphi=l\varphi N$ implies $\lambda=l\sigma(\lambda)$, which is
only possible if $\lambda=0$. Hence $N=0$ on $E$. Finally, we
necessarily have $\varphi(1,0)=(1,0)$ and $\varphi(0,1)=(\alpha,1)$
for some $\alpha\in K_0$ and conversely giving $\alpha\in K_0$
uniquely determines~$\varphi$.

Denote by $K_0^v$ the $\QQl$-vector space underlying~$K_0$. To
$\alpha\in K_0^v$ associate the unique extension~(\ref{eq:sesfpm})
with $E=K_0\oplus K_0$ and $\varphi(0,1)=(\alpha,1)$. This determines
a surjective map
$$K_0^v\map{\longrightarrow}{j}\Ext^1_{\MFNf{K}}(K_0,K_0)\;,$$ and one
checks that it is in fact $\QQl$-linear. Take $x\in K_0$ and let $L$
be the automorphism of $K_0\oplus K_0$ given by~$\bigl({1\ x\atop 0\
1}\bigr)$. If we equip the source with $\varphi$ and the target with
$\varphi'$, then $L$ induces an equivalence of the associated
extensions if $\varphi'=L^{-1}\varphi L$. It follows that
$j(\alpha)=j(\alpha+\sigma(x)-x)$, so the kernel of $j$ is a
sub-$\QQl$-vector~space of $K_0^v$ of codimension~1; this implies that
$\dim_{\QQl}\Ext^1_{\MFNf{K}}(K_0,K_0)\le1$.

But on the other hand,
$\dim_{\QQl}\Ext^1_{\unrrep{G_K}}(\QQl,\QQl)=1$. To see this, suppose
$V$ is unramified and sits in an extension of $\QQl$ by $\QQl$.
Taking a suitable basis for $V$, the action of the Galois group is
given by $\bigl({1\ \eta\atop\hphantom{\eta}\ 1}\bigr)$, where $\eta$
is an unramified character $\Gal(\QQlbar/K)\to\QQl$. In other words,
$\eta$ is a morphism of groups $\ZZhat\to\QQl$ which, being
continuous, must factor through $\ZZl$. But $\Hom(\ZZl,\QQl)$ is
one-dimensional.

Thus we obtain that $V$ in~(\ref{eq:sesrepr}) is unramified when
restricted to $\Gal(\QQlbar/K)$. Then if $g$ is an element of the
inertia subgroup of~$\Gal(\QQlbar/\QQl)$, the $[K:\QQl]$-th power of
$g$ must act trivially. But on the other hand $g$ must act
unipotently, as $V$ sits in the short exact
sequence~(\ref{eq:sesrepr}). Hence $g$ itself must act trivially and
$V$ is an unramified representation of $\Gal(\QQlbar/\QQl)$.
\end{proof}

Now let
$$0=H_0\subset H_1\subset \cdots\subset H_{P_i}=H$$ be a
Jordan-H\"older filtration of~$H$. Hence $H_j/H_{j-1}\simeq\QQl$ and
each $H_j$ is unramified outside $l$ and potentially semi-stable
at~$l$.  Clearly, $H_1\simeq\QQl$; assume that $H_j\simeq\QQln{j}$ for
some $j$~~($1\le j\le P_i-1$). Then $H_{j+1}$ is everywhere unramified
by the above lemma. However, as Minkowski's theorem says that $\QQ$
has no non-trivial unramified extensions, we conclude that $H_{j+1}$
is isomorphic to~$\QQln{j+1}$. So by induction,
$H\simeq\QQln{P_i}$. This finishes the proof of Thm.~\ref{thm:main}.

\section{The Hodge structure}

Let $V$ be an integral, regular, quasi-projective
$\QQ$-scheme. Consider a finite group~$G$ acting on~$V$ from the
right. Denote by $f\colon V\to \Xpc$ the canonical projection to the
quotient scheme~$\Xpc=V/G$. Note that there is a natural $G$-action on
the module~$f_*\Omega^1_{V/\QQ}$.

We are interested in the cohomology of the quotient space~$\Xpc$ and
in particular in the Hodge structure of the singular cohomology of its
associated analytic space. For this we will use~\cite{Steenbrink}. In
order to be able to apply some results of~\cite{Steenbrink} we will
need the following result:

\begin{proposition}\label{prop:diff}
Let $\Sigma\subset \Xpc$ be a closed subset of codimension~$\ge2$
containing all singular points of~$\Xpc$. Let $j\colon
U\hookrightarrow \Xpc$ be the open subscheme complementary
to~$\Sigma$. Then, for all~$p$, there is a canonical isomorphism
\begin{equation}\label{eq:Gequiv}
j_*\Omega^p_{U/\QQ}\map\longrightarrow\sim (f_*\Omega_{V/\QQ}^p)^G
\end{equation}
of $\shfO_\Xpc$-modules.
\end{proposition}

\begin{proof}
For any $\QQ$-scheme~$Z$, we abbreviate $\Omega_Z:=\Omega_{Z/\QQ}^p$.
Put $W=f^{-1}(U)$ and let $g\colon W\to U$ be the restriction of~$f$.
Note that $G$ acts on $W$ with quotient~$U$.  The canonical map
$g^*\Omega_{U}\to\Omega_{W}$ induces a morphism
$$\alpha\colon \Omega_{U}\map\longrightarrow\sim(g_*g^*\Omega_{U})^G
  \longrightarrow (g_*\Omega_{W})^G\;.$$ This map is an isomorphism in
  the stalk at the generic point.  To define (\ref{eq:Gequiv}), we
  apply $j_*$ to $\alpha$ to obtain a map from $j_*\Omega_U$ to the
  sheaf $j_*(f_*\Omega_V)^G$. This last sheaf is isomorphic to
  $(f_*\Omega_V)^G$, as a consequence of the fact that if $Z$ is a
  regular $\QQ$-variety and $z\colon Z'\hookrightarrow Z$ is a dense
  open subset with complement of codimension~$\ge2$, then
  $z_*\Omega_{Z'}\simeq\Omega_Z$. Another consequence of this fact is
  that the map~(\ref{eq:Gequiv}) thus defined is independent of the
  choice of~$\Sigma$. In particular, we have reduced the problem to
  showing that for any point $\eta\in U$ whose closure has
  codimension~$1$, the map $\alpha$ is an isomorphism in the stalk
  at~$\eta$.

The question being local for the \'etale topology on~$U$ and $V$, it
suffices to consider, for $d\ge1$ and $n\ge0$, the action of $\MU_d$
on $\Affine^n_\QQ$, with the quotient map
$\Affine_\QQ^n\to\Affine_\QQ^n$ mapping $(a_1,\ldots,a_n)$ to
$(a_1^d,a_2,\ldots,a_n)$. The results follows from an easy
calculation.
\end{proof}

Let $j\colon U\hookrightarrow \Xpc$ be the smooth locus of~$\Xpc$.  As
a consequence of Prop.~\ref{prop:diff}, there is a canonical
isomorphism between the De~Rham complexes
$$j_*\Omega^\bullet_{U/\QQ}\map\longrightarrow\sim
  (f_*\Omega^\bullet_{V/\QQ})^G$$ Therefore, for each~$i$ we obtain an
  isomorphism
\begin{equation}\label{eq:DRisom}
\rH^i(\Xpc,j_*\Omega^\bullet_{U/\QQ})\map\longrightarrow\sim
  \rH^i(\Xpc,(f_*\Omega^\bullet_{V/\QQ})^G)\map\longrightarrow\sim
  \rH^i_\DR(V/\QQ)^G\;.
\end{equation}
This last isomorphism follows from the fact that $f$ is finite and
that taking cohomology with $\QQ$-coefficients commutes with taking
invariants under a finite group action. The Hodge filtration on the
De~Rham complex induces filtrations on the vector spaces
in~(\ref{eq:DRisom}) and the isomorphisms in~(\ref{eq:DRisom}) respect
those filtrations.

\medskip

Fix a prime~$l$. Denote by $\Ddr$ the Fontaine functor (see
e.g.,~\cite{Tsuji}) from the category of
$\Gal(\QQlbar/\QQl)$-representations over~$\QQl$ to the category of
finite dimensional filtered $\QQl$-vector spaces.

\begin{proposition}\label{prop:comparison}
In the above situation, suppose furthermore that $\Xpc$ is proper
over~$\QQ$.  For every~$i$, there is an isomorphism of filtered vector
spaces
$$\Ddr(\rH^i(\Xpc_{\QQlbar,\et},\QQl))\map\longrightarrow\sim
  \rH^i(\Xpc,j_*\Omega^\bullet_{U/\QQ})\otimes_\QQ\QQl\;.$$
\end{proposition}

\begin{proof}
As $V$ is proper and smooth, the comparison theorem (see
e.g.,~\cite[Thm.~A1]{Tsuji}) states that
$\Ddr(\rH^i(V_{\QQlbar,\et},\QQl))$ and $\rH^i_\DR(V/\QQl)$ are
isomorphic. This isomorphism preserves the $G$-invariants.  Using
(\ref{eq:DRisom}) we obtain the isomorphism of the proposition.
\end{proof}

With these prerequisites, we are finally ready to prove the corollary
to the main theorem:

\begin{corollary}\label{cor:main}
In the situation of Theorem~\ref{thm:main} (with $U$ a non-empty open
subscheme of $\Spec\ZZ$), suppose that the coarse moduli space
$\Xcoarse$ of $\X_\QQ$ is the quotient of a smooth projective
$\QQ$-scheme by a finite group\footnote{This last condition can be
omitted, by invoking \cite[Cor.~8.12]{Bogaart} directly after the
statement, in the proof, that $\Ddr(W)\simeq\Ddr(W^{\rm ss})$ if
$W$~is de~Rham. }.

Then for each~$i$, there is an isomorphism of $\QQ$-Hodge structures
$$\rH^i(\Xcoarse(\CC),\QQ)\simeq \cases{0,&if $i$ is odd,\cr
  \QQ(-i/2)^{P_{i/2}},&if $i$ is even,}$$ where the left hand side is
  equipped with the canonical Hodge structure of~\cite{Hodge3}.
\end{corollary}

\begin{proof}
Consider a short exact sequence
  $$0\to W'\to W\to W''\to 0$$
in the category of representations of $\Gal(\QQlbar/\QQl)$ over~$\QQl$.
If $W$~is a de~Rham representation, then the corresponding sequence
  $$0\to\Ddr(W')\to\Ddr(W)\to\Ddr(W'')\to0$$
is a short exact sequence on the underlying vector spaces and all maps
are strict with respect to the filtrations (see \cite[3.4,~3.7
and~3.8]{FontaineAst}). 
It follows that $\Ddr(W)\simeq\Ddr(W^{\rm ss})$ if $W$~is de~Rham.

Without loss of generality we assume that $\Xcoarse$ is integral.  Note that
$\rH^i(\X_{\QQlbar,\et},\QQl)$ is isomorphic to
$\rH^i(\Xcoarse_{\QQlbar,\et},\QQl)$ by Lemma~\ref{lemma:cohomcoarse}
and note that these representations are de~Rham.
Therefore, combining Thm.~\ref{thm:main} and
Prop.~\ref{prop:comparison} it suffices to exhibit an isomorphism of
filtered vector spaces
$$\rH^i(\Xcoarse(\CC),\CC)\map\longrightarrow\sim
\rH^i(\Xcoarse,j_*\Omega^\bullet_{U/\QQ}) \otimes_\QQ\CC\;,$$ where
$j\colon U\to\Xcoarse$ is the smooth locus of~$\Xcoarse$.  This is
done in \cite[Thm.~1.12]{Steenbrink}.
\end{proof}

\end{document}